\documentclass[12pt]{amsart}

\usepackage{latexsym, amsmath, amscd, amssymb, amsthm} 
\usepackage[T2A]{fontenc}
 \usepackage[cp1251]{inputenc}
\newtheorem{te}{Теорема}
\newtheorem{pro}{Пропозиція}
 \newtheorem{oz}{Означення}

\begin{document}

\noindent
УДК 512.815

 \title{Аналоги прямо\"i та обернено\"i теореми Робертса для тернарних форм }

\author{Л. Бедратюк}\address{ Хмельницький національний університет, вул. Інститутська, 11, 29016, м. Хмельницький}

\maketitle

{\small 

 Для тернарних    форм   доводяться  аналоги добре відомої  в теорії інваріантів теореми Робертса. Встановлено, що незвідні  коваріанти, контраваріанти та змішані конкомітанти тернарної форми однозначно визначаються їхніми  старшими  членами.
}

{\small 
 Analogues of invariant theory's well-known Roberts theorem are proved for ternary forms. We established that covariants, contravariants and mixed concomitants of a ternary form are uniquely determined by their lead coefficients.
}

\section{Вступ}

Розглянемо $\mathbb{K}$-векторний простір $T_n$ тернарних \\ форм степеня $n:$
$$ 
u(x_1,x_2,x_3)=\sum_{i+j\leq n} \, \frac{n!}{i! j! (n{-}(i+j))!}a_{i, j} x_1^{n-(i+j)} x_2^i x_3^j,
$$
де $a_{i, j} \in \mathbb{K}$, a $\mathbb{K}$ -- поле  нульової характеристики.
Координатне кільце $R_n$ простору  $T_n$ ототожнимо з алгеброю многочленів $k[A]:=\mathbb{K}[a_{0,\, 0},a_{1,\, 0},\ldots, a_{0,\, n}]$ від $\displaystyle \frac{1}{2} (n+1)(n+2)$ змінних, а координатне кільце  простору $T_n\oplus \mathbb{K}^3\oplus (\mathbb{K}^3)^*$ ототожнимо з кільцем многочленів  $\mathbb{K}[A, X, U]:=\mathbb{K}[A, x_1, x_2, x_3, u_1, u_2,u_3].$ 
Стандартна дія групи  $SL_3$ підстановками на $T_n$  індукує дію  $SL_3$ і на  кільці $\mathbb{K}[A,X,U]$.
 Поліноміальні функції з $R_n,$ та  $\mathbb{K}[A,X,U],$ які залишаються інваріантними відносно дії групи $SL_3$ утворюють кільця $R_n^{SL_3},$ та $\mathbb{K}[A,X,U]^{SL_3},$ які називаються, відповідно, кільцями інваріантів  та змішаних конкомітантів  тернарної форми степеня $n.$ Кільце $\mathbb{K}[A,X]^{SL_3}$ та кільце $\mathbb{K}[A,U]^{SL_3}$ називаються  кільцями коваріантів та  контраваріантів тернарної форми степеня $n$ ( див. \cite{GrYo}).  Зокрема,  форма $u(x_1,x_2,x_3)$ буде коваріантом степеня $n.$ Для довільного многочлена, який є однорідним по кожному наборі змінних $A,$  $X$ та $U,$  його степені відносно цих наборів   називаються  відповідно степенем, порядком  та класом.

 Знаходження явного вигляду комітантів --  породжуючих елементів вищеозначених кілець інваріантів,   є основною  задачею класичної теорії інваріантів, яка була розв'язана  ще  Горданом \cite{Gor-1}, але лише для $n\leq 3.$ Найвищим досягненням того періоду було  обчислення в докторській дисертації Е.~Ньотер \cite{Net}  мінімальної системи із 331 породжуючих   кільця інваріантів $\mathbb{K}[A]^{SL_3}$ для $n=4.$
 Майже всі відомі  конкомітанти отримані в неявному вигляді  символічним методом, коли конкомітанти  зображуються через  трансвектанти, тобто як результат дії деякого $SL_3$-інваріантного диференціального оператора   на  конкомітанти   менших степенів. 

Одним із підходів до вивчення комітантів могло би бути   встановлення аналогу теореми Робертса для тернарних форм. В класичному формулюванні теорема Робертса стверджує (див. \cite{Gle}, \cite{SP}), що всякий коваріант бінарної форми степеня $n$   відносно дії гру-пи $SL_2$ однозначно визначається своїм старшим членом -- коєфіцієнтом біля $x_1^n.$ В свою чергу, старший член всякого коваріанта  бінарної  форми  є  інваріантом     одновимірної підалгебри верхніх трикутних матриць,  іншими словами, він є старшим вектором деякого незвідного $sl_2$-модуля.  Тому проблема опису кільця коваріантів бінарних  форм зводиться до питання опису  кільця інваріантів підалгебри верхніх трикутних  матриць в алгебрі Лі $sl_2.$

В даній роботі для комітантів тернарних форм степеня $n$  доведені твердження, які є аналогами теореми Робертса. Показано, що коєфіцієнти  незвідного  комітанта  породжують незвідний  $sl_3$-модуль  в $k[A]$ і старший коєфіцієнт комітанта буде старшим вектором цього модуля.   Справедливе і обернене твердження -- всякий інваріант алгебри  верхніх трикутних матриць $UT_3$ є старшим вектором  деякого  $sl_3$-модуля. Таким чином, встановлені твердження зводять задачу   знаходження породжуючих елементів кільця  $\mathbb{K}[A,X,U]^{SL_3}$ до  простішої задачі знаходження породжуючих елементів кільця $\mathbb{K}[A]^{UT_3}.$ 
 
\section{Елементи Казиміра }

Дамо означення   елемента Казиміра, який  буде головним обчислювальним засобом при вивченні комітантів тернарної форми.
\begin{oz}
Симетричним добутком $U \cdot V$ векторних просторів $V$ та $U$ назвемо підалгебру    симетричної алгебри $S(U \oplus V),$   породжену елементами вигляду $u\, v,$  $ v \in V,$ $u \in U$
\end{oz}

Якщо простори $U,$ $V$ є $sl_3$-модулями то їхній симетричний добуток  $U \cdot V$ також буде  $sl_3$-модулем, якщо покласти $$g(uv)=g(u)\,v+u g(v),$$  для всіх $g \in sl_3,$ $ v \in V,$ $u \in U.$

\begin{oz}
Всякий інваріант   $sl_3$-модуля $U \cdot V$ називається елементом Казиміра.
\end{oz}

\begin{te}\label{null}
Припустимо, що $U,$ $V$ -- два  $sl_3$-модулі. В $sl_3$-модулі $U \cdot V$ елемент Казиміра існує  тоді і лише тоді, коли $U\cong V^*.$
\end{te}
\begin{proof}[Доведення] Припустимо,що  $U\cong V^*,$  $m=\dim U.$ Виберемо в просторах $V$ і $V^*$ дуальні базиси $\{v_i\},$ $\{v^*_i\},$ ${ i=1\ldots m.}$ Довільний  елемент $z \in sl_3$ діє як лінійний оператор в  $V$ і $V^*$.  Добре відомо, що матриці $C{=\{ c_{i\,j} \},}$ ${C^*=\{ c^*_{i\,j} \}}$ цього оператора в дуальних  базисах зв'язані співвідношенням $C^*=(-C)^{\rm T}$. Покажемо, що елемент $$v_1 v^*_1+v_2 v^*_2+\cdots +v_m v^*_m \in V \cdot V^*,$$   є інваріантом. Маємо
$$
\begin{array}{l}
\displaystyle z\Bigl(\sum_{i=1}^{m} v_i v_i^*\Bigr)=\sum_{i=1}^{m} (z(v_i) v_i^*+v_i\,z(v_i^*))=
\displaystyle=\sum_{i=1}^{m} \Bigl(\sum_{j=1}^{m}c_{i j}v_j v_i^*+v_i z(v_i^*)\Bigr)=
\\
\displaystyle =\sum_{i=1}^{m} \Bigl(\sum_{j=1}^{m}c_{j i}v_i v_j^*+v_i\,z(v_i^*)\Bigr)=
\displaystyle=\sum_{i=1}^{m}  v_i\Bigl(\sum_{j=1}^{m}c_{i j}v_j^*+z(v_i^*)\Bigr)=0. \\
\end{array}
$$

Припустимо тепер, що елемент $$v_1 u_1+v_2 u_2+\cdots +v_m u_m \in V \cdot U$$  є інваріантом. Аналогічно знаходимо  
$$
\begin{array}{l}
\displaystyle z\Bigl(\sum_{i=1}^{m} v_i u_i\Bigr)=
\sum_{i=1}^{m}  v_i\Bigl(\sum_{j=1}^{m}c_{j i}u_j+z(u_i)\Bigr). \\
\end{array}
$$
Рівність нулю можлива лише тоді,  коли для всіх $i$ буде виконуватися  $\displaystyle  \sum_{j=1}^{m}c_{j i}u_j+z(u_i)=0,$ тобто дія $z$ на $U$ є контрагредієнтною до дії на $V,$  а це означає, що   $U\cong V^*.$
\end{proof}
Елемент Казиміра $sl_3$-модуля $U \cdot V$  будемо позначати $\Delta(U,V).$  Якщо в просторах $U,$ $V$  задано контрагредієнтні базиси  $$U:=\langle u_1, u_2,\ldots u_m \rangle,  \mbox{     }  V:=\langle v_1, v_2,\ldots v_m \rangle,$$ то 
$$
\Delta(U,V)=u_1 v_1+u_2 v_2+\cdots +u_m v_m.
$$
Можна показати, що елемент Казиміра не залежить від вибору пар дуальних базисів в просторах $U$ i $V^*.$

\section{Реалізація $sl_3$-модулів в $k[A].$}

Незвідний  $sl_3$-модуль  з старшою вагою $[m_1,m_2]$ будемо позначати $\Gamma_{m_1,\,m_2},$ або, бажаючи явно вказати старший вектор $v$ -- $ \Gamma_{m_1,\,m_2}(v).$ Розмірність  простору $\Gamma_{m_1,\,m_2}$  рівна  (див. \cite{FH}) $$ \frac{1}{2} (m_1+1) (m_2+1) (m_1+m_2+2).$$
Старший вектор $sl_3$-модуля $\Gamma_{m_1,\,m_2}$ є інваріантом підалгебри $UT_3$ верхніх трикутних матриць. Аналогічно молодший вектор цього модуля є інваріантом підалгебри $DT_3$ нижніх трикутних матриць. Якщо вектор  $u$ є  старшим вектором  старшої ваги $[m_1,m_2],$ то справедливе співвідношення  (див. \cite{Hum}) $$\Gamma_{m_1,m_2}(u)=\mathfrak{U}(DT_3)(u),$$  де через $\mathfrak{U}(L)$ позначено  універсальну  огортуючу алгебру алгебри Лі $L.$ Аналогічно
 для молодшого вектора  $v$   ваги  $[-m_1,-m_2]$  отримаємо $${\Gamma_{m_1,m_2}(v)=\mathfrak{U}(UT_3)(v)}.$$ 
 Позначимо через $E_{i\,j},$ $i,j=1 \ldots 3,$ матричні одиниці, тобто такі матриці у яких на перетині $i$-го рядка та $j$-го стовпчика знаходиться одиниця, а на всіх інших місцях нулі. Має місце співвідношення
 $$[E_{i\,j},E_{k\,l}]:=E_{i\,j}E_{k\,l} -E_{k\,l}E_{i\,j}=\delta_{j\,k} E_{i\,l}-\delta_{i\,l} E_{k\,j}.$$ 
Матриці $E_{1\,2},$ $E_{2\,3},$ $E_{1\,3}$ утворюють базис підалгебри $UT_3$ верхніх трикутних матриць, матриці $E_{2\,1},$ $E_{3\, 2},$ $E_{3\, 1}$  утворюють базис підалгебри $DT_3$  нижніх трикутних матриць. Матриці $E_{1\,1}{-}E_{2\, 2},$ $E_{2\, 2}{-}E_{3\,3},$  та $E_{1\,1}{-}E_{3\, 3},$ породжують картанівську підалгебру   в $sl_3.$ 

Для довільного $sl_3$-модуля $V$ позначимо через $D_1,$ $D_2,$ $D_3$ лінійні оператори з ${\rm End}(V),$ які відповідають дії на $V$ відповідно елементів $E_{1\,2},$ $E_{2\,3},$ $E_{1\,3}.$ Аналогічно оператори $\hat D_1,$ $\hat D_2,$ $\hat D_3$ відповідають дії на $V$ відповідно елементів $E_{2\,1},$ $E_{3\, 2},$ $E_{3\, 1},$  оператори $E_1,$ $E_2,$ i $E_3$  відповідають дії  відповідно елементів $E_{1\,1}{-}E_{2\, 2},$ $E_{2\, 2}{-}E_{3\,3},$ $E_{1\,1}{-}E_{3\, 3}.$

Випишемо комутаційні співвідношення між цими операторами, які  нам будуть потрібні  в подальшому:
$$
\begin{array}{lll}
[E_1,D_1]=2 D_1,  & [E_1,D_3]=D_3, & D_2, \hat D_1]=0,\\

 [E_1,\hat D_1]=-2 \hat D_1,  & [E_1, \hat D_3]=-\hat D_3, & [E_1,\hat D_2] = \hat D_2,\\

[D_1, \hat D_1]=E_1,  & [D_1, \hat D_3]=- \hat D_2, & [D_1, \hat D_2] =0,\\

[E_1,D_2] =-D_2,[  & [D_2, \hat D_3]=- \hat D_1,  & [D_2, \hat D_2] =E_2, \\

[D_3, \hat D_1]=- D_2,    & [D_3, \hat D_3]=-D_1, & [D_3, \hat D_2] =E_2. \\

\end{array}
$$

Алгебра $sl_3$ діє на векторному просторі $$X:=\langle x_1, x_2, x_3 \rangle, $$ диференціюваннями, а саме 
$$
\begin{array}{ll}
\displaystyle D_1=-x_2 \frac{\partial}{\partial x_1}, &  \displaystyle D_2=-x_3 \frac{\partial}{\partial x_2}, \\
E_1=\displaystyle  x_2 \frac{\partial}{\partial x_2}-x_1\frac{\partial}{\partial x_1}, &E_2= \displaystyle  x_3 \frac{\partial}{\partial x_3}-x_2\frac{\partial}{\partial x_2},   \\
\displaystyle \hat D_1=-x_1 \frac{\partial}{\partial x_2}, &  \displaystyle \hat D_2=-x_2 \frac{\partial}{\partial x_3}, \\
 \displaystyle D_3=-x_3 \frac{\partial}{\partial x_1}, & E_3= \displaystyle  x_3 \frac{\partial}{\partial x_3}-x_1\frac{\partial}{\partial x_1},
\end{array}
$$
$$
\begin{array}{lll}
 \displaystyle &  \displaystyle \hat D_3=-x_1 \frac{\partial}{\partial x_3}.  &
\end{array}
$$
Векторний простір $X$   є стандартним незвідним $sl_3$-модулем ізоморфним до $\Gamma_{0,\,1},$ а векторний простір $U:=\langle u_1, u_2, u_3 \rangle \cong X^* $  є  незвідним $sl_3$-модулем ізоморфним до $\Gamma_{1,\,0}.$ Відповідний елемент Казиміра 
$$
u:=\Delta(X,U)=x_1 u_1+x_2 u_2+x_3 u_3,
$$
називається універсальним коваріантом.

Симетричні степені $S^m(X)$ і  $S^m(U)$  є незвідними $sl_3$-модулями ізоморфними відповідно до  $\Gamma_{0\!,m}$ і $\Gamma_{m\!,0}.$
	
Вивчимо  дію алгебри $sl_3$ на породжуючі елементи кільця $R_n.$
\begin{pro}\label{pro1}
В $sl_3$-модулі $R_n$ відповідні диференціальні оператори діють за формулами
$$
\begin{array}{ll}
\displaystyle \! \! D_1(a_{i,j})=i\,a_{i{-}1,j}, & \! \! D_2(a_{i,j})=j\,a_{i{+}1,j{-}1},\\
\displaystyle \! \! \hat D_1(a_{i,j})=(n-(i+j))\,a_{i{+}1,j}, & \! \!  \hat D_2(a_{i,j})=i\,a_{i{-}1,j{+}1},\\
  \hat \! \! D_3(a_{i,j})=(n-(i+j)) a_{i,j+1}, & \! \! D_3(a_{i,j})=j a_{i,j-1}, \\
\! \! E_1(a_{i,j})=(n-(2i+j)) a_{i,j}, & \! \! E_2(a_{i,j})=(i-j) a_{i,j},
\end{array}
$$
$$E_3(a_{i,j})=(d-(i+2j)) a_{i,j}. $$
\end{pro}
\begin{proof}[Доведення]
Для доведення використаємо той факт, що форма $u(x_1,x_2, x_3)$ є коваріантом і тому кожен з операторів  $D_i, 
\hat D_i$ повинен зануляти її. 
Зокрема, для диференціюваня $D_1,$  маємо 
$$
D_1(u(x_1,x_2,x_3))=\sum_{i+j\leq n} \, \frac{n!}{i! j! (n{-}(i+j))!}(D_1(a_{i,j}) x_1^{n-(i+j)} x_2^i x_3^j+a_{i,j} D_1(x_1^{n-(i+j)} x_2^i x_3^j))=
$$
$$
=D_1(a_{0,1}) x_1^{n-1} x_3+\cdots+D_1(a_{0,n}) \frac{1}{n!} x_3^n+ \sum_{\begin{array}{c} \mbox{\small \it i+j}\leq n \\  i>0 \end{array} } \Bigl( D_1(a_{i,j})-i\,a_{i{-}1,j}\Bigr) x_1^{d-(i+j)} x_2^i x_3^j.
$$
Тому  рівність $D_1(u(x_1,x_2,x_3))=0$  можлива лише за умови, що всі коєфіцієнти  рівні нулю, отже, маємо $D_1(a_{0,j})=0$ для всіх $0\leq j \leq n,$ і   $D_1(a_{i,j})=i\,a_{i{-}1,j} ,$ що і потрібно було показати.  В такий самий спосіб визначається дія  інших операторів на кільці $R_n.$
\end{proof}

Якщо елемент $a \in K[A]$ є власним вектором   оператора $E_i,$ $i=1,2,$ то його власне значення будемо познaчати $\omega_i(a)$ і називати $i$-вагою елемента $a,$ а такий елемент $a$  будемо називати ваговим вектором.  Зрозуміло, що $i$-вага є лінійною, адитивною функцією на множині вагових векторів.  Однорідний многочлен $a$ з $K[A]$ називається ізобарним, якщо він буде ваговим відносно обох операторів $E_1,$ $E_2.$ У цьому випадку набір $[\omega_1(a), \omega_2(a)]$ буде називатися вагою многочле\-на $a.$
\begin{te}\label{null1}
Нехай  $V:=\{ v_k\}$  $V^*:=\{ v^*_k\}$ -- два дуальні $sl_3$-модулі, причому всі базисні вектори є ваговими векторами. Якщо для деякого номера $i$  елемент  $v_i \in V$  є старшим вектором  у $V,$ то $v^*_i$ буде молодшим вектором у $V^*.$ 
\end{te}
\begin{proof}[Доведення]

Оскільки базисні вектори є ваговими, то ваговим буде  і елемент $v_i v^*_i,$   його вага рівна $$[\omega_1(v_i)+\omega_1(v^*_i),\omega_2(v_i)+\omega_2(v^*_i)].$$ Із умов $E_1(\Delta(V,V^*))=E_2(\Delta(V,V^*))=0$ випливає що і $E_1(v_i v^*_i)=E_2(v_i v^*_i)=0$ звідки знаходимо, що  $\omega_1(v^*_i)=-\omega_1(v_i),$  $\omega_2(v^*_i)=-\omega_2(v_i).$ Нехай для деякого номера $i$ елемент $v_i$ є  старшим вектором  $sl_3$-модуля $V$ старшої ваги $[\omega_1(v_i),\omega_2(v_i)].$ Тоді вектор $v^*_i$ має вагу  $-[\omega_1(v_i),\omega_2(v_i)]$  і ця вага буде молодшою  вагою $sl_3$-модуля $V^*.$
\end{proof}

Наступні теореми встановлюють правила обчислень в $sl_3$-модулі $\mathfrak{U}(UT_3)(a),$ де $a$ є старшим вектором в $k[A].$
\begin{pro}\label{pro2}  Нехай $a$ -- однорідний, ізобарний многочлен  з $k[A].$ Тоді
$$
E_1(\hat D_1^{\alpha} \hat D_2^{\beta} \hat D_3^{\gamma} (a))=(\omega_1(a)-2\alpha+\beta -\gamma) \hat D_1^{\alpha} \hat D_2^{\beta} \hat D_3^{\gamma} (a),
$$
$$
E_2(\hat D_1^{\alpha} \hat D_2^{\beta} \hat D_3^{\gamma} (a))=(\omega_2(a)+\alpha-2 \beta +\gamma) \hat D_1^{\alpha} \hat D_2^{\beta} \hat D_3^{\gamma} (a).
$$
\end{pro}
\begin{proof}[Доведення] $(i)$ Використовуючи комутаційне співвідношення $[E_1,\hat D_1]=-2\hat D_1$ отримаємо, що 
$$E_1 \hat D_1 (a)=[E_1,\hat D_1](a)+D_1(E_1(a))=$$
$$=-2 \hat D_1(a)+\omega_1(a) \hat D_1(a)=(\omega_1(a)-2) \hat D_1(a).$$
В загальному випадку   маємо, що $$E_1 \hat D_1^{\alpha} (a)=(\omega_1(a)-2\alpha) \hat D_1^{\alpha}(a).$$ 

Враховуючи співвідношення  $[E_1,\hat D_2]=\hat D_2,$ та \\ $[E_1,\hat D_3]=- \hat D_3$   знаходимо, що $$E_1 \hat D_2^{\beta} (a)=(\omega_1(a)+\beta) \hat D_2^{\beta}(a),$$
$$E_1 \hat D_3^{\gamma} (a)=(\omega_1(a)-\gamma) \hat D_3^{\gamma}(a).$$
В загальному випадку отримаємо 
$$
E_1(\hat D_1^{\alpha} \hat D_2^{\beta} \hat D_3^{\gamma} (a))=(\omega_1(a)-2\alpha+\beta -\gamma) \hat D_1^{\alpha} \hat D_2^{\beta} \hat D_3^{\gamma} (a).
$$
\noindent
$(ii)$ Аналогічними   міркуваннями, враховуючи співвідношення 
$$[E_2,\hat D_1]= D_1, [E_2,\hat D_2]=-2\hat D_2, [E_2,\hat D_3]=- \hat D_3, $$ знаходимо
$$
E_2(\hat D_1^{\alpha} \hat D_2^{\beta} \hat D_3^{\gamma} (a))=(\omega_2(a)+\alpha-2\beta -\gamma) \hat D_1^{\alpha} \hat D_2^{\beta} \hat D_3^{\gamma} (a).
$$
\end{proof}
Як наслідок отримуємо, що довільний незвідний $sl_3$-модуль $V$ в $k[A]$ із старшим вектором $a$ старшої  ваги $[d_1,d_2]$ розкладaється в суму  вагових підпросторів $V_{(i,\,j)}$, де $$V_{(i,\,j)}=\{\hat D_1^{\alpha} \hat D_2^{\beta} \hat D_3^{\gamma} (a), \lambda_1=i, \lambda_2=j  \}.$$
Тут $\lambda_1:=\omega_1(a)-2\alpha+\beta-\gamma, $ $\lambda_2:=\omega_2(a)+\alpha-2\beta -\gamma$
i $ [-d_1,-d_2] \leq [i,j] \leq [d_1,d_2].$
\begin{pro}\label{pro3}  Нехай $a$ -- однорідний ізобарний многочлен з $k[A]^{UT_3}.$ Тоді
$$
\begin{array}{l}
\displaystyle D_1(\hat D_1^{\alpha} \hat D_2^{\beta} \hat D_3^{\gamma} (a)){=}\alpha (\lambda_1{+}\alpha+1) \hat D_1^{\alpha-1} \hat D_2^{\beta} \hat D_3^{\gamma} (a)-\gamma \hat D_1^{\alpha} \hat D_2^{\beta+1} \hat D_3^{\gamma-1} (a), \\
\displaystyle D_2(\hat D_1^{\alpha} \hat D_2^{\beta} \hat D_3^{\gamma} (a))=\beta (\omega_2(a)- \beta +1) \hat D_1^{\alpha} \hat D_2^{\beta-1} \hat D_3^{\gamma} (a)+\gamma \hat D_1^{\alpha+1} \hat D_2^{\beta} \hat D_3^{\gamma-1} (a). 
\end{array}
$$
\end{pro}
\begin{proof}[Доведення]
Доведемо лише першу формулу, оскільки друга формула доводиться за тією ж схемою. Маємо $$D_1 \hat D_3(a)=[D_1, \hat D_3](a)+\hat D_3 D_1(a)=-\hat D_2(a).$$
За індукцією неважко показати, що $${D_1 \hat D_3^{\gamma}(a)=-\gamma \hat D_2 \hat D_3^{\gamma-1}(a),}$$  і, враховуючи комутативність операторів $D_1$ i $\hat D_2,$   знаходимо $${D_1 \hat D_2^{\beta} \hat D_3^{\gamma}(a)=-\gamma \hat D_2^{\beta+1} \hat D_3^{\gamma-1}(a)}.$$

Для довільного однорідного ізобарного $a'$ маємо $$D_1 \hat D_1(a'){=}[D_1, \hat D_1](a')+\hat D_1 D_1(a'){=}\omega_1(a')+\hat D_1 D_1(a').$$  В загальному випадку  отримаємо
$$
\begin{array}{l}
\displaystyle D_1 \hat D_1^{\alpha}(a')=\displaystyle(\sum_{\tau=0}^{\alpha-1}\omega_1(\hat D_1^{\tau}(a')) D_1^{\alpha-1}(a')+D_1^{\alpha} D_1(a')=\displaystyle \alpha (\omega_1(a')-\alpha+1) D_1^{\alpha-1}(a')+\hat D_1^{\alpha} D_1(a').
\end{array}
$$
Підставивши $a'= \hat D_2^{\beta} \hat D_3^{\gamma}(a),$ отримаємо необхідне співвідношення.
\end{proof}

\section{Коваріанти тернарної форми}

Наступне твердження є аналогом відомої теореми Робертса про коваріанти бінарної форми.
\begin{te}
Нехай  
$$
 f{=}\sum_{i+j\leq d} \, \frac{d!}{i! j! (d{-}(i+j))!}b_{i,j} x_1^{d-(i+j)} x_2^i x_3^j,  \mbox{  } b_{i\,j} \in k[A]
$$
-- незвідний коваріант порядку $d$.
Тоді :
\begin{enumerate}
\item[({\it i})] векторний простір $B_d:=\langle \{b_{i,j} \}, i+j\leq d \rangle $ буде незвідним $sl_3$-модулем ізоморфним до         $\Gamma_{d,\,0}.$ 

\item[({\it ii})] елемент $b_{0,\,0}$ буде старшим вектором $sl_3$-модуля $B_d$ з старшою вагою $[d,0].$

\item[({\it iii})] коваріант $f$ є елементом Казимира, \\${f=\Delta(B_d,S^d(X))}$ i записується  у вигляді
$$
f=\sum_{i+j\leq d} \, \frac{1}{i! j!}\hat D_1^i \hat D_3 ^j( b_{0,0})  x_1^{d-(i+j)} x_2^i x_3^j.
$$
\end{enumerate}
\end{te}
\begin{proof}[Доведення] $(i)$  Очевидно, що $f$ є інваріантом $sl_3$-модуля $ B_d \cdot S^d(X).$  Тому із теореми \ref{null} отримуємо $B_d \cong S^d(X)^* \cong \Gamma_{d,\,0}.$ 

\noindent
$(ii)$
Аналогічно, як і у пропозиції \ref{pro1} знаходимо дію на $B_d$ диференціальних операторів, які відповідають породжуючим елементам алгебри $sl_3:$ 
$$
\begin{array}{ll}
\displaystyle D_1(b_{i,j})=i\,b_{i{-}1,j}, & D_2(b_{i,j})=j\,b_{i{+}1,j{-}1},\\
\displaystyle \hat D_1(b_{i,j})=(d-(i+j))\,b_{i{+}1,j}, & \hat D_2(b_{i\,j})=i\,b_{i{-}1,j{+}1}.
\end{array}
$$
Очевидно, що в $B_d$ є лише один інваріант алгебри $DT_3,$  а саме $b_{0,0},$  тому $b_{0,0}$  є старшим вектором $sl_3$-модуля $B_d$ з вагою $[d,0].$

\noindent
$(iii)$  Оскільки  $b_{0\,0}$  є старшим вектором незвідного модуля $B_d$, то $B_d=\mathfrak{U}(DT_3)(b_{0\,0})$  ( див., наприклад, \cite{Hum}). Тому всі $b_{i\,j}$ можна виразити  через степені операторів $\hat D_1,$ $\hat D_2,$  i $\hat D_2.$ Використовуючи явний вигляд дії цих операторів неважко показати, що 
$$
b_{i,j}=\frac{1}{m (m-1) \ldots (m-(i+j-1)} \hat D_1^i \hat D_3^j (b_{0,\,0}).
$$

\noindent
Базиси $\displaystyle \left \{\frac{d!}{i! j! (d{-}(i+j))!}\,b_{i,j} \right  \}{=}\left \{\frac{1}{i! j!}\hat D_1^i \hat D_3 ^j( b_{i,j}) \right  \}$ та \\$ \left \{x_1^{d-(i+j)} x_2^i \,x_3^j \right \},$ $ i+j \leq d$ є взаємно дуальними, тому 
$$
f{=}\Delta(B_d,S^d(X)){=}\sum \, \frac{1}{i! j!}\hat D_1^i \hat D_3 ^j( b_{0,0}) x_1^{d-(i+j)} x_2^i\, x_3^j.
$$
Тут сума береться по таких індексах $i$, $j,$ для яких $i+j\leq d.$
\end{proof}

Отже,  всякий незвідний коваріант однозначно визначається своїм старшим коєфіцієнтом, який є інваріантом підалгебри $UT_3,$ утвореної верхньо-три-кутними матрицями.

 Для формулювання оберненої теореми введемо поняття порядку многочлена відносно операторів $\hat D_1,$ $\hat D_2,$.
\begin{oz}
Для довільного многочлена \\ $z \in \mathbb{K}[A,X,U]$ набір $[z]:=[{\rm ord}_1(z),{\rm ord}_2(z)],$  де $${\rm ord}_i(z):=\max\{s, \hat D_i^s(a) \neq 0  \}, i=1,2$$ називається порядком $z$ відносно  оператора $\hat D_i.$
\end{oz}
Цілі числа  ${\rm ord}_i(z),$ $i=1,2$ будемо називати $i$-порядком многочлена $z$ відносно диференціювання $\hat D_i.$
Коректність означення порядку випливає із того, що оператори диференціювання $\hat D_1,$ $\hat D_2$  є локально нільпотентними в кільці $\mathbb{K}[A,X,U].$ 

Покажемо, що для однорідних, ізобарних елементів з $\mathbb{K}[A]^{UT_3}$  їхні порядки співпадають з відповідними вагами. Справедливе  наступне  твердження
\begin{pro}
Нехай $a$ однорідний, ізобарний  елемент з  $\mathbb{K}[A]^{UT_3}.$ Тоді $[a]=[\omega_1(a),\omega_2(a)].$
\end{pro}
\begin{proof}[Доведення] Покажемо, що $\omega_1(a)={\rm ord}_1(a).$ За означенням порядку елемента маємо, що 
$$ \hat D_1^{\displaystyle{\rm ord}_1(a)+1}(a)=0.$$  З іншого боку, враховуючи пропозицію \ref{pro2} отримаємо 
$$
\begin{array}{l}
\displaystyle D_1 ( \hat D_1^{{\rm ord}_1(a){+}1}(a))=\\
\displaystyle=({\rm ord}_1(a)+1) (\omega_1(a)-{\rm ord}_1(a)) \hat D_1^{{\rm ord}_1(a)}(a)=0.
\end{array}
$$
Отже, $\omega_1(a)={\rm ord}_1(a),$  що і потрібно було довести. 
Аналогічно показується, що $\omega_2(a)={\rm ord}_2(a).$
\end{proof}

\begin{te}
Нехай $a$ -- однорідний, незвідний, ізобарний  елемент з  $\mathbb{K}[A]^{DT_3}$ порядку $[d,0].$ Тоді
\begin{enumerate}
\item[({\it i})] векторний простір $$\overline B_d:={\mathfrak U}(UT_3) a=\displaystyle \left \{\frac{1}{[d,i+j-1]!}\hat D_1^i \hat D_3 ^j( a), i+j\leq d \right \},$$ є незвідним $sl_3$-модулем ізоморфним до $\Gamma_{d,\,0}.$
\item[({\it ii})] елемент Казиміра $\Delta(\bar B_d,S^d(X))$ є коваріантом порядку $d$ тернарної форми.
\end{enumerate}
\end{te}
\begin{proof}[Доведення] $(i)$ Покладемо $$\bar b_{i,j}=\displaystyle \frac{1}{[d,i+j-1]!}\hat D_1^i \hat D_3 ^j( a).$$  Тоді
$$
\begin{array}{l}
\displaystyle D_1(\bar b_{i,j})=D_1\left (\frac{1}{[d,i+j-1]!}\hat D_1^i \hat D_3 ^j( a)\right )=\\
\displaystyle =\frac{i (\omega_1(a)-i-j+1)}{[d,i+j-1]!}  \hat D_1^{i-1} \hat D_3 ^j( a)-j \hat D_1^i \hat D_2 \hat D_3 ^{j-1}( a).
\end{array}
$$
Враховуючи  комутативність   операторів $\hat D_2$ та $\hat D_3$ і те, що ${\rm ord}_2(a)=0,$ тобто $\hat D_2(a)=0,$  отримаємо, що другий доданок рівний нулю. Взявши до уваги  те, що $\omega_1(a)=d,$ отримаємо
$$
\begin{array}{l}
\displaystyle D_1(\bar b_{i,j})=\frac{i}{[d,i+j-1]!} (d-i-j+1) \hat D_1^{i-1} \hat D_3 ^j( a)=\\
\displaystyle =\frac{i}{[d,i+j-2]!}  \hat D_1^{i-1} \hat D_3 ^j( a)=i \bar b_{i-1,j}.
\end{array}
$$
Далі, 
$$
\begin{array}{l}
\displaystyle \hat D_1(\bar b_{i,j})=\frac{1}{[d,i+j-1]!}\hat D_1^{i+1} \hat D_3 ^j( a)=\\
\displaystyle =\frac{d-(i+j)}{[d,i+j]!}\hat D_1^{i+1} \hat D_3 ^j( a)=(d-(i+j)) \bar b_{i+1,j}.
\end{array}
$$
Аналогічно  знаходимо, що $D_2(b_{i,j})=j\, \bar b_{i+1,j-1}$ i \\ $\hat D_2(\bar b_{i,j})=i\, \bar b_{i-1,j+1}.$  Таким   чином, векторний простір $\bar B_d$ є $sl_3$-модулем, причому дія алгебри $sl_3$  співпадає з її дією   на $sl_3$-модулі $R_d,$ див. пропозицію \ref{pro1}.  Оскільки розмірності просторів $\bar B_d$ i $R_d$ рівні $\displaystyle \frac{1}{2} (d+1)(d+2),$  то $\bar B_p \cong R_d \cong \Gamma_{d,0}.$

\noindent
$(ii)$ Розглянемо  білінійну форму  $$(\cdot,\cdot ): \bar B_d \times S^d(X) \to \mathbb{K},$$  значення якої на базисних елементах відповідних просторів визначається за формулою 
$$
(\bar b_{k\,l}, x_1^{d-(i+j)} x_2^i x_3^j)=\frac{i! j! (d-(i+j))!}{d!} \delta_{i\,k} \delta_{j\,l},
$$
тут $\delta_{i\,j}$ -- символ Кронекера.

Перевіримо, що ця форма є $sl_3$-інваріантною,  тоб-то для всіх $g\in sl_3,$  $u \in \bar B_d,$ $v \in S^d(X)$
 має місце співвідношення 
$
(g(u),v)+(u,g(v))=0.
$

Для оператора $D_1$ маємо
$$
\begin{array}{l}
\bigl (D_1(\bar b_{k\,l}),x_1^{d-(i+j)} x_2^i x_3^j \bigr )=
\displaystyle k (\bar b_{k-1\,l},x_1^{d-(i+j)} x_2^i x_3^j)= \\
\displaystyle{=}k \frac{i! j! (d-(i+j))!}{d!} \delta_{k-1\,i} \delta_{l\,j}{=}
\displaystyle(i+1) \frac{i! j! (d{-}(i+j))!}{d!},
\end{array}
$$
i
$$
\begin{array}{l}
(\bar b_{k\,l},D_1(x_1^{d-(i+j)} x_2^i x_3^j))
\displaystyle{=}{-} \frac{(\bar b_{k\,l},x_1^{d{-}(i+j+1)} x_2^{i+1} x_3^j)}{(d{-}(i+j))^{-1}}=\\
\displaystyle{=}{-}(d{-}(i{+}j)) \frac{(i{+}1)! j! (d{-}(i+j+1))!}{d!} \delta_{k\,i+1} \delta_{l\,j}=\\
\displaystyle=-\frac{(i+1)! j! (d-(i+j))!}{d!}.
\end{array}
$$

Таким чином,
$$
\bigl (D_1(\bar b_{k\,l}),x_1^{d-(i+j)} x_2^i x_3^j\bigr)+\bigl (\bar b_{k\,l},D_1(x_1^{d-(i+j)} x_2^i x_3^j)\bigr)=0.
$$
Аналогічно виконується перевірка інваріантності і для операторів $\hat D_1,$ $D_2$ i $\hat D_2.$ Отже, білінійна форма $(\cdot,\cdot )$ є $sl_3$-інваріантною і тому базиси 
$$
\left \{\frac{d!}{i! j! (d-(i+j))!} \bar b_{i,j} \right \} \mbox{   та    } \bigl \{x_1^{d-(i+j)} x_2^i x_3^j \bigr\} , i+j \leq d,
$$
взаємно дуальні.

 Відповідний елемент Казиміра є $sl_3$-інваріантом  i  має вигляд 
$$
\Delta(\bar B_d, S^d(X)){=}\sum_{i+j\leq d} \, \frac{n!}{i! j! (d{-}(i+j))!}\bar b_{i,j} x_1^{d-(i+j)} x_2^i x_3^j,
$$
тобто є коваріантом степеня $d.$
\end{proof}

\section{Контраваріанти тернарної форми}

Перейдемо до вивчення контраваріантів тернарної форми. Аналогічно, як і у випадку коваріантів мають місце наступні теореми
\begin{te}
Нехай  
$$
 f=\sum_{i+j\leq d} \, \frac{d!}{i! j! (d{-}(i+j))!}c_{i\,j} u_3^{d-(i+j)} u_1^i u_2^j, \mbox{  } c_{i,j} \in k[A], 
$$
-- незвідний контраваріант порядку $d$.
Тоді 
\begin{enumerate}
\item[({\it i})] векторний простір $C_d:=\langle \{c_{i,j} \}, i+j\leq d \rangle $ є незвідним $sl_3$-модулем ізоморфним $\Gamma_{0,\,d}.$ 

\item[({\it ii})] елемент $c_{0,0}$ є старшим вектором $sl_3$-моду-\\ля $C_d$ з вагою $[0,d].$

\item[({\it iii})] контраваріант $f$ є елементом Казимира  $${f=\Delta(C_d,S^d(U)),}$$ який записується  у вигляді
$$
f=\sum_{i+j\leq d} \, \frac{(-1)^{i+j}}{i! j!}\hat D_2^j \hat D_3 ^i( c_{0,0})  u_3^{n-(i+j)} u_1^i u_2^j.
$$
\end{enumerate}
\end{te}
\begin{proof}[Доведення]
$(i)$  Очевидно, що $f$ є інваріантом $sl_3$-модуля $ C_d \cdot S^d(U).$  Тому, із теореми \ref{null} отримуємо $$C_d \cong S^d(X)^* \cong \Gamma_{0,\,d}.$$

\noindent
$(ii)$
Аналогічно, як і у пропозиції \ref{pro1} знаходимо дію на $C_d$ диференціальних операторів, які відповідають породжуючим елементам алгебри $sl_3:$ 
$$
\begin{array}{ll}
\displaystyle \! \! \! D_1(c_{i,j}){=}-i\,с_{i{-}1,j+1}, & \! \! D_2(c_{i,j}){=}-j\,c_{i,j{-}1},  \\
\displaystyle \! \! \! \hat D_2(c_{i,j}){=}-(d-(i+j))\,c_{i,j{+}1}, &\! \! \hat D_1(c_{i,j}){=} {-}j\,c_{i{+}1,j{-}1},\\
\displaystyle  \! \! \! \hat D_3(c_{i,j}){=}-(d-(i+j))\,c_{i+1,j}, & \! \! D_3(c_{i,j}){=}-i\,c_{i-1,j}.
\end{array}
$$
Очевидно, що в $C_d$ є лише один інваріант алгебри $DT_3,$  а саме $c_{0,0},$  тому $c_{0,0}$  є старшим вектором незвідного  $sl_3$-модуля $C_d$ з вагою $[0,d].$

\noindent
$(iii)$  Оскільки елемент $c_{0,0}$  є старшим вектором незвідного $sl_3$-модуля $C_d$ , то $C_d=\mathfrak{U}(DT_3)(c_{0,0}).$   Тому всі $c_{i,j}$ можна виразити  через степені операторів $\hat D_1,$ $\hat D_2,$  i $\hat D_3.$ Використовуючи явний вигляд дії цих операторів неважко показати, що 
$$
b_{i,j}=\frac{(-1)^{i+j}}{[m,i+j-1]!} \hat D_2^j \hat D_3^i (c_{0,0}).
$$
\noindent
Базиси $\left \{\displaystyle \frac{d! \,c_{i,j}}{i! j! (d{-}(i+j))!} \right \}{=}\left \{\displaystyle \frac{(-1)^{i+j}}{i! j!}\hat D_1^i \hat D_3 ^j( c_{0,0}) \right \}$ та $\bigl \{u_3^{d-(i+j)} u_1^i \, u_2^j \bigr \},$ $ i+j \leq d$ є взаємно дуальними, тому 
$$
\Delta(B_d,S^d(U)){=}\sum_{i+j\leq d} \, \frac{(-1)^{i+j}\hat D_1^i \hat D_3 ^j(c_{0,0}) u_3^{d-(i+j)} u_1^i \, u_2^j}{i! j!}.
$$

\end{proof}
Отже, кожен незвідний контраваріант  однозначно вичначає своїм старшим членoм, який є інваріантом алгебри $DT_3.$ В наступній теоремі  доводиться справедливість оберненого твердження.
\begin{te}
Нехай $a$ є незвідний, однорідний, ізобарний  елемент з  $k[A]^{UT_3}$ мультипорядку $[0,d].$ То-ді
\begin{enumerate}
\item[({\it i})] векторний простір  $$\overline C_d:={\mathfrak U}(UT_3) a=\displaystyle \left \{\frac{(-1)^{i+j}}{[d,i+j-1]!}\hat D_2^j \hat D_3 ^i( a), i+j\leq d  \right \},$$ є незвідним $sl_3$-модулем ізоморфним до $\Gamma_{0,\,d}.$
\item[({\it ii})] елемент Казиміра $\Delta(\bar C_d,S^d(U))$ є контраваріантом порядку $d$ тернарної форми.
\end{enumerate}
\end{te}
\begin{proof}[Доведення]\! \! $(i)$ Покладемо $\bar c_{i,j}{=}\! \displaystyle \frac{(-1)^{i+j}\hat D_2^j \hat D_3 ^i(a)}{[d,i+j-1]!}.$  То-ді, використавши пропозицію \ref{pro2}, знайдемо

$$
\begin{array}{l}
\displaystyle D_1(\bar c_{i,j})=D_1\left(\frac{(-1)^{i+j}}{[d,i+j-1]!}\hat D_2^j \hat D_3 ^i( a)\right)=\\
\displaystyle =\frac{-(-1)^{i+j} i}{[d,i+j-1]!}\hat D_2^{j+1} \hat D_3 ^{i-1}( a)=-i \bar c_{i-1,j+1}.
\end{array}
$$

Неважко переконатися, використовуючи індукцію, що $\hat D_1 \hat D_2^j(a)=-j \hat D_3 \hat D_2^{j-1}.$  Тому, знову взявши до уваги пропозицію \ref{pro2}, враховуючи комутативність операторів $\hat D_1,$ i $\hat D_3,$ а, також те, що мають місце рівності $${\omega_2(a){=}{\rm ord}_2(a){=}d},$$  після нескладних обчислень отримаємо
$$
\begin{array}{l}
\displaystyle D_2(\bar c_{i,j})=D_2\left(\frac{(-1)^{i+j} \hat D_2^j \hat D_3 ^i( a)}{[d,i+j-1]!}\right)= \\
\displaystyle =\frac{ -(-1)^{i+j-1}j }{{[d,i+j-2]!} }  \hat D_2^{j-1}\hat D_3 ^{i}( a)) =-j \bar c_{i,j-1}.
\end{array}
$$
Далі, аналогічно знаходимо  
$$
\hat D_1(\bar c_{i,j})=\frac{-j (-1)^{i+j} \hat D_2^{j-1} \hat D_3 ^{i+1}( a)}{[d,i+j-1]!}=-j \bar c_{i+1,j-1},
$$
і
$$
\hat D_2(\bar c_{i,j})=-(d-(i+j)) \bar c_{i,j+1}.
$$
Отже, векторний простір $\bar C_d$ є $sl_3$-модулем, причому дія алгебри $sl_3$  співпадає з її дією   на $sl_3$-модулі $C_d,$ див.  $(ii)$  попередньої теореми.  Оскільки розмірності просторів $\bar C_d$ i $C_d$ рівні $\displaystyle \frac{1}{2} (d+1)(d+2),$  то $\bar C_d \cong C_d \cong \Gamma_{0,d}.$

\noindent
$(ii)$ Розглянемо  білінійну форму  $$(\cdot,\cdot ): \bar C_d \times S^d(U) \to \mathbb{K},$$  значення якої на базисних елементах відповідних просторів визначається за формулою 
$$
\bigl (\bar c_{k,l}, u_3^{d-(i+j)} u_1^i u_2^j\bigr )=\frac{(-1)^{i+j} i! j! (d-(i+j))!}{d!} \delta_{i\,k} \delta_{j\,l}.
$$
Аналогічно, як і у випадку коваріантів можна показати, що ця форма є $sl_3$-інваріантною і тому базиси 
$$
\left \{\frac{(-1)^{i+j} d!}{i! j! (d-(i+j))!} \bar c_{i,j} \right \} \mbox{   та    } \bigl \{u_3^{d-(i+j)} u_1^i u_2^j) \bigr \} , i+j \leq d,
$$
будуть взаємно дуальними.

 Відповідний елемент Казиміра є $sl_3$-інваріантом  i  має вигляд 
$$
\Delta(\bar C_d, S^d(U))=\sum_{i+j\leq n} \, \frac{(-1)^{i+j} d!}{i! j! (d{-}(i+j))!}\bar b_{i,j} u_3^{d-(i+j)} u_1^i u_2^j,
$$
тобто є контраваріантом порядку $d.$
\end{proof}

\section{Змішані конкомітанти тернарної форми}

Перейдемо до вивчення змішаних конкомітантів тернарної форми. Як і у випадку коваріантів та контраваріантів  мають місце наступні теореми
\begin{te}
Нехай  
$$
 f=\sum_{\begin{array}{c} \mbox{{\small i+j}} \leq d_1 \\  \mbox{{\small k+l}} \leq d_2 \end{array} } \, \frac{ d_1! d_2! B_{i\!,j}^{k\!,l} x_1^{d_1-(i+j)} x_2^i x_3^j u_3^{d_2-(k+l)} u_1^k u_2^l}{i! j! k! l! (d_1{-}(i+j))! (d_2-(k+l))!} , 
$$
-- незвідний змішаний конкомітант  класу  $[d_1,d_2].$ Тут $B_{i\!,j}^{k\!,l} \in k[A]$  i $d_1, d_2 >0.$ 
Тоді:
\begin{enumerate}
\item[({\it i})] векторний простір $$B_{d_1}^{d_2}:=\langle \{B_{i\!,j}^{k\!,l} \}, i+j\leq d_1, k+l\leq d_2 \rangle $$ є незвідним $sl_3$-модулем ізоморфним  $\Gamma_{d_1,d_2}.$ 

\item[({\it ii})] елемент $B_{0\!, 0}^{0\!, 0}$ є старшим вектором $sl_3$- модуля $B_{d_1}^{d_2}$ з вагою $[d_1,d_2].$

\item[({\it iii})] змішаний конкомітант $f$ є елементом Казимира, причому  $f=\Delta(B_{d_1}^{d_2},S^{d_1}(X) \cdot S^{d_2}(U)).$ 
\end{enumerate}
\end{te}
\begin{proof}[Доведення]
$(i)$  Очевидно, що $f$ є інваріантом $sl_3$-модуля $ B_{d_1}^{d_2} \cdot S^{d_1}(X) \cdot S^{d_2}(U).$  Тому із теореми \ref{null} отримуємо $$ B_{d_1}^{d_2}  \cong (S^{d_1}(X) \cdot S^{d_2}(U))^* =(S^{d_1}(X))^* \cdot (S^{d_2}(U))^*.$$
$sl_3$-модуль $S^{d_1}(X) \cdot S^{d_2}(U)$ можна легко розкласти на незвідні підмодулі. Для цього знайдемо старші вектори в $S^{d_1}(X) \cdot S^{d_2}(U),$ тобто вектори інваріантні відносно алгебри $DT_3$.  Неважко переконатися, що за старші вектори можна взяти такі многочлени
$$v_i:=x_3^{d_1-i} u_1^{d_2-i} u^i, \mbox{   } i \leq 0\ldots i_0, \mbox{ }  i_0:={\min}(d_1,d_2), $$
 a $u:=x\,u_1+y\,u_2+z\,u_3$ -- універсальний коваріант.

Оскільки ${\rm ord}_1(u_1){=}1,$ і  ${{\rm ord}_1(x_3){=}{\rm ord}_1(u)=0,}$  то \\${\rm ord}_1(v_i)=d_1-i.$ Аналогічно знаходимо, що i    \\${{\rm ord}_2(v_i)=d_2-i.}$ Оскільки,  порядок інваріанта спі-впадає з його вагою, то  вага кожного вектора $v_i$ рів-на ${[d_1{-}i,d_2{-}i],}$  і  має місце розклад
$$
S^{d_1}(X) \cdot S^{d_2}(U)=\Gamma_{d_1,\,d_2} (v_0)+ \cdots +\Gamma_{d_1-i_0,\,d_2-i_0}(v_{i_0}).
$$
Тут $\Gamma_{d_1\!{-}i,\,d_2\!{-}i}(v_i)$ -- незвідний $sl_3$-модуль із старшим вектором $v_i$ і вагою ${[d_1{-}i,\, d_2{-}i]}.$ 

Повністю аналогічно можна перевірити, що многочлени $\bar v_i:=x^{d_1-i} u_3^{d_2-i} u^i$ є молодшими векторами $sl_3$-модуля $S^{d_1}(X) \cdot S^{d_2}(U).$

Задамо лінійне відображення  векторних просторів $\varphi :S^{d_1}(X) \cdot S^{d_2}(U) \to B_{d_1}^{d_2},$  яке кожному елементу з $S^{d_1}(X) \cdot S^{d_2}(U)$ ставить у відповідність дуальний йому елемент з $B_{d_1}^{d_2},$  тобто 
$$
\varphi\bigl ( x_1^{d_1-(i+j)} x_2^i x_3^j u_3^{d_2-(k+l)} u_1^k u_2^l\bigr)=\frac{ d_1! d_2!}{i! j! k! l! (d_1{-}(i+j))! (d_2-(k+l))!} B_{i\,j}^{k\,l}.
$$
Врахувавши теорему \ref{null1}, отримаємо, що відображеня $\varphi$  переводить молодший вектор $sl_3$-модуля $S^{d_1}(X) \cdot S^{d_2}(U)$  у старший вектор $sl_3$-модуля $B_{d_1}^{d_2}.$ Позначимо через $w_i:=\varphi(v_i)$ відповідні старші вектори, а через $\Gamma_{d_1-i,d_2-i}(w_i)$ -- відповідні незвідні підмодулі. Має місце розклад
$$
B_{d_1}^{d_2}=\Gamma_{d_1,d_2} (w_0)+ \cdots +\Gamma_{d_1-i_0,d_2-i_0}(w_{i_0}).
$$
Розглянемо тепер   тепер елемент Казиміра $$\Delta(B_{d_1}^{d_2},S^{d_1}(X) \cdot S^{d_2}(U)).$$  Очевидно має місце розклад 
$$
\Delta(B_{d_1}^{d_2},S^{d_1}(X) \cdot S^{d_2}(U)){=}
\Delta(\Gamma_{d_1,d_2} (w_0),\Gamma_{d_1,d_2} (v_0))+ \cdots +\Delta(\Gamma_{d_1-i_0,d_2-i_0}(v_{i_0}),\Gamma_{d_1-i_0,d_2-i_0}(w_{i_0})) .
$$
Кожен із елементів Казиміра в правій частині, крім останього, є змішаним конкомітантом, тому  елемент $\Delta(B_{d_1}^{d_2},S^{d_1}(X) \cdot S^{d_2}(U))$ буде незвідним лише тоді, коли всі змішані конкомітанти правої частини меншого класу будуть рівні нулю, а це можливо лише у випадку, коли всі старші вектори $w_i=0,$ $i=1\ldots i_0,$ а $w_0 \neq 0.$ Таким чином отримаємо
$$
\Delta(B_{d_1}^{d_2},S^{d_1}(X) \cdot S^{d_2}(U))=\Delta(\Gamma_{d_1,d_2} (w_0),\Gamma_{d_1,d_2} (v_0)),
$$
звідки зразу випливає, що $B_{d_1}^{d_2}=\Gamma_{d_1,d_2} (w_0) \cong \Gamma_{d_1,d_2}.$

\noindent
$(ii)$

Аналогічно до того, як це було зроблено в пропозиції \ref{pro1},  знаходимо дію  породжуючих елементів алгебри $sl_3$ на $B_{d_1}^{d_2}:$
$$
\begin{array}{l}
 \hat D_1(B_{i,j}^{k,l}){=}(d_1-(i+j)) B_{i+1,j}^{k,l}{-}l B_{i,j}^{k+1,l-1} ,\\
       \hat D_2(B_{i,j}^{k,l}){=}i B_{i-1,j+1}^{k,l}{-}(d_2-(k+l)) B_{i,j}^{k,l+1} , \\
\hat D_3(B_{i,j}^{k,l}){=}(d_1-(i+j)) B_{i,j+1}^{k,l}{-}(d_2-(k+l)) B_{i,j}^{{k+1},l},\\
D_1(B_{i,j}^{k,l})=i B_{i-1,j}^{k,l}-k B_{i,j}^{k-1,l+1} ,\\
D_2(B_{i,j}^{k,l})=j B_{i+1,j-1}^{k,l}-l B_{i,j}^{k,l-1} ,\\
D_3(B_{i,j}^{k,l})=j B_{i,j-1}^{k,l}-k B_{i,j}^{k-1,l}.
\end{array}
$$
Звідси зразу отримуємо, що $$D_1(B_{0,0}^{0,0})=D_2(B_{0,0}^{0,0})=0,$$  отже, $B_{0,0}^{0,0}$ -- старший вектор у $B_{d_1}^{d_2}.$

\noindent
$(iii)$  Оскільки базиси $$\displaystyle \left  \{\frac{ d_1! d_2!}{i! j! k! l! (d_1{-}(i+j))! (d_2-(k+l))!} B_{i,j}^{k,l}\right \} \mbox{ та } \bigl \{x_1^{d_1-(i+j)} x_2^i x_3^j u_3^{d_2-(k+l)} u_1^k u_2^l \bigr\}, i+j \leqslant d_1, k+l \leqslant d_1$$ дуальні, то $f=\Delta(B_{d_1}^{d_2},S^{d_1}(X) \cdot S^{d_2}(U)).$  
\end{proof}
Нехай $a$ є незвідний, однорідний, ізобарний  елемент з  $k[A]^{UT_3}$ порядку $[d_1,d_2].$ Розглянемо векторний простір $\bar B_{d_1}^{d_2}$ породжений елементами  $\left \{\bar B_{i,j}^{k,l} \right \},$  де $\bar B_{0,0}^{0,0}{:=}a,$  а всі $\bar B_{i\,j}^{k\,l}$ знаходяться із умов
$$
\begin{array}{l}
 \hat D_1(\bar B_{i,j}^{k,l})=(d_1-(i+j)) \bar B_{i+1,j}^{k,l}-l \bar B_{i,j}^{k+1,l-1} ,\\
          \hat D_2(\bar B_{i,j}^{k,l})=i \bar B_{i-1,j+1}^{k,l}-(d_2-(k+l)) \bar B_{i,j}^{k,l+1} , \\
  \hat D_3(\bar B_{i,j}^{k,l}){=}(d_1-(i+j)) \bar B_{i,j+1}^{k,l}{-}(d_2-(k+l)) \bar B_{i,j}^{k+1,l},
\end{array}
$$
та  із умов рівності нулю елементів $\bar w_i,$ $i=1\ldots i_0.$
 Тут $\bar w_i$ утворюються із $w_i$ заміною $B_{i,j}^{k,l}$ на $\bar B_{i,j}^{k,l}.$

\begin{te}
Нехай $a$ є незвідний, однорідний, ізобарний  елемент з  $k[A]^{UT_3}$ мультипорядку $[d_1,\, d_2].$ Тоді:
\begin{enumerate}
\item[({\it i})] векторний простір  $\bar B_{d_1}^{d_2}:={\mathfrak U}(UT_3) a$ є незвідним $sl_3$-модулем ізоморфним до $\Gamma_{d_1,\,d_2}.$
\item[({\it ii})] елемент Казиміра $\Delta(\bar B_{d_1}^{d_2},S^{d_1}(X) \cdot S^{d_2}(U))$   є змішаним конкомітантом тернарної форми класу $[d_1,d_2].$
\end{enumerate}
\end{te}
\begin{proof}[Доведення] $(i)$ Пряма перевірка показує, що $\bar B_{d_1}^{d_2}$ є $sl_3$-модулем. Оскільки простір $\bar B_{d_1}^{d_2}$ побудований таким чином, що в ньому існує   лише один старший вектор ваги $[d_1,d_2]$, а саме $\bar B_{0,0}^{0,0}=a,$ то $\bar B_{d_1}^{d_2}$ є незвідним $sl_3$-модулем, причому $\bar B_{d_1}^{d_2}={\mathfrak U}(UT_3) a.$

\noindent
$(ii)$ Розглянемо  білінійну форму  $$(\cdot,\cdot ): \bar B_{d_1}^{d_2} \times S^{d_1}(X) \cdot S^{d_2}(U) \to \mathbb{K},$$  значення якої на базисних елементах відповідних просторів визначається за формулою 

$$
\left (\bar B_{i\,j}^{k\,l}, x_1^{d_1-(i+j)} x_2^{i'} x_3^{j'} u_3^{d_2-(k'+l')} u_1^{k'} u_2^{l'}\right )=\frac{i! j! k! l! (d_1{-}(i+j))! (d_2-(k+l))!}{ d_1! d_2!} \delta_{i,i'}  \delta_{j,j'}  \delta_{k,k'}  \delta_{k,k'}.
$$
Аналогічно, як і у випадку коваріантів  та контраваріантів можна показати  що ця форма є невиродженою та $sl_3$-інваріантною. Тому, базиси 
$$
\left \{\frac{ d_1! d_2!}{i! j! k! l! (d_1{-}(i+j))! (d_2-(k+l))!} \bar B_{i\,j}^{k\,l}  \right \},
$$
та 
$$
\{x_1^{d_1-(i+j)} x_2^{i} x_3^{j} u_3^{d_2-(k+l)} u_1^{k} u_2^{l} \} , i+j \leq d_1, k+l \leq d_2,
$$
є дуальними, а 
відповідний елемент Казиміра $$\Delta(\bar B_{d_1}^{d_2},S^{d_1}(X) \cdot S^{d_2}(U)),$$  є змішаним конкомітантом тернарної форми класу $[d_1,d_2].$
\end{proof}

{\bf Приклад.}  Для  $n=3$,   розглянемо в $k[A]$ многочлен
$
a:=a_{0,\,0} a_{2,\,0}-a_{1,\,0}^2
$
який є інваріантом підалгебри $k[A]^{DT_3}.$
Оскільки $\hat D_1^3(a)=0,$ але $\hat D_1^2(a) \neq0,$ то ${\rm ord}_1(a)=2.$ Аналогічними міркуваннями   отримуємо, що ${\rm ord}_2(a)=2,$  отже, 
 порядок $a$ рівний $[2,2],$  і $a$ буде  старшим вектором зі старшою вагою $[2,2]$ $sl_3$-підмодуля $\bar B_2^2={\mathfrak U}(UT_3) a$ в $k[A],$ ізоморфного стандартному $sl_3$-модулю $\Gamma_{2,\,2}.$
Відповідно до вагової діаграми знайдемо базисні вектори вагових підпросторів $B_{(i,j)}$  $sl_3$-модуля ${\mathfrak U}(UT_3) a$
$$
\begin{array}{ll}
B_{(2,2)}=\{B_{0\,0}^{0\,0}=a\}, & \phantom{**********}  \\ B_{(0,3)} =\{ \hat D_1(a)=2 B_{1 \,0}^{0\,0} \}, & \\  B_{(3,0)}=\{D_2(B_{0\,0}^{0\,0})=-2 B_{0\,0}^{0\,1}\}, & \\
B_{(4,-2)} =\{\hat D_2^2(a) = 2 B_{0,0}^{0,2} \}, &   \\ B_{(-2,4)} =\{ \hat D_1^2(a)=2 B_{2 \,0}^{0\,0} \}, &
\\   B_{(-3,3)} =\{  \hat D_1^2 \hat D_3(a) = -4 B_{2,0}^{1,0} \}, & \\
B_{(-4, 2)}=\{ \hat D_1^2 \hat D_3^2(a) = 4 B_{2,0}^{2,0} \}, & \\ B_{(3,-3)}=\{ \hat D_2^2 \hat D_2(a) = 4 B_{0,1}^{0,2} \}, &\\  B_{(2,-4)}=\{\hat D_2^2 \hat D_3^2(a)= \{ 4 B_{0,2}^{0,2} \} ,& \\
 B_{(-2,-2)}=\{ \hat D_3^4(a) = 24 B_{0,2}^{2,0} \}, & \\
B_{(1,1)} =\{ \hat D_3(a) =  2 B_{0,1}^{0,0}-2 B_{0,0}^{1,0},  \\ \hat D_1 D_2(a) = -4 B_{1,0}^{0,1}+2 B_{0,0}^{1,0} \}, &\\
B_{(-1, 2)}=\{ \hat D_1 \hat D_3(a) = -4 B_{1,0}^{1,0}+2 B_{1,1}^{0,0}, &\\   \hat D_1^2 \hat D_2 = 8 B_{1,0}^{1,0}-4 B_{2,0}^{0,1} \} , &\\
B_{(2, -1)} =\{\hat D_1 \hat D_2^2(a) = 4B_{1,0}^{0,2}-4 B_{0,0}^{1,1}, & \\ \hat D_2 \hat D_3(a) = 2 B_{0,0}^{1,1}-4 B_{0,1}^{0,1} \}, & \\
B_{(-2, 1)} =\{\hat D_1 \hat D_3(a) =4 B_{1,0}^{2,0}-8 B_{1,1}^{1,0}, & \\ \hat D_1^2 \hat D_2  \hat D_3(a)=-8 B_{1,0}^{2,0}+4 B_{2,0}^{1,1}+8 B_{1,1}^{1,0} \}, & \\
B_{(1,-2)} =\{\hat D_2 \hat D_3^2(a)=8 B_{0,1}^{1,1}-4 B_{0,2}^{0,1}, & \\ \hat D_1 \hat D_2^2 \hat D_3(a) = 4 B_{1,1}^{0,2}-8 B_{0,1}^{1,1}\}, & \\
B_{(-1,-1)} =\{ \hat D_1 \hat D_2 \hat D_3^2(a)=8 B_{1,1}^{1,1}-8 B_{0,1}^{2,0}+4 B_{0,2}^{1,0}, &\\  \hat D_1^2 \hat D_2^2 \hat D_3(a) = -16 B_{1,1}^{1,1}+8 B_{0,1}^{2,0}  \}, & \\
B_{(-3,0)}=\{ \hat D_1^2 \hat D_2 \hat D_3^2(a) = -16 B_{1,1}^{2,0} \}, & \\ B_{(0,-3)}=\{ \hat D_1 \hat D_2^2 \hat D_3^2(a) = -8 B_{0,2}^{1,1} \}, & \\
B_{(0,0)} =\{ \hat D_3^2(a)=-8 B_{0,1}^{1,0}+2 B_{0,0}^{2,0}+2 B_{0,2}^{0,0}, & \\  \hat D_1 \hat D_2 \hat D_3(a)=4 B_{1,0}^{1,1}-2 B_{0,0}^{2,0}-4 B_{1,1}^{0,1}+4 B_{0,1}^{1,0}, &\\
 \hat D_1^2 \hat D_2^2(a)=-16 B_{1,0}^{1,1}+4 B_{0,0}^{2,0}+4 B_{2,0}^{0,2} \}.
\end{array}
$$

Отримали $27$ рівнянь для  $36$  невідомих $B_{i,j}^{k,l}.$ Інші $9$ рівнянь знайдемо із наступних  міркувань.

В $sl_3$-модулі $S^{2}(X) \cdot S^{2}(U)$ розглянемо молодші вектори 
$$
\begin{array}{l}
v_1=x_1 u_1 u=x_1 u_1(x_1 u_1+x_2 u_2+x_3 u_3)=x_1^2 u_1^2+x_1 x_2 u_1 u_2+x_1 x_3 u_1 u_3,\\

v_2=u^2=(x_1 u_1+x_2 u_2+x_3 u_3)^2.
\end{array}
$$ 
Тоді вектори 
 $$
\begin{array}{l}
\varphi(v_1)=B_{0,1}^{0,0}+B_{0,0}^{1,0}+B_{1,0}^{0,1}, \\
\varphi(v_2)=B_{0,0}^{2,0}+B_{2,0}^{0,2}+B_{0,2}^{0,0}+2 B_{1,0}^{1,1}+2 B_{0,1}^{1,0}+2 B_{1,1}^{0,1},
\end{array}
$$ 
будуть старшими векторами $B_2^2$  мультипорядків $[1,1]$ i $[0,0].$ Розмірність $\Gamma_{1,1}(\varphi(v_1))$ рівна $8$ і розмірність $\Gamma_{0,0}(\varphi(v_2))$ рівна $1.$ 
Одномірні вагові підпростори $\Gamma_{1,1}(\varphi(v_1))$ породжуються такими елементами

$$
\begin{array}{l}
\varphi(v_1)= B_{0,1}^{0,0}+B_{0,0}^{1,0}+B_{1,0}^{0,1},  \\  \hat D_1 (\varphi(v_1))= B_{1,0}^{1,0}+B_{2,0}^{0,1}+B_{1,1}^{0,0}, \\
\hat D_1 \hat D_3(\varphi(v_1)) = -B_{1,0}^{2,0}-B_{2,0}^{1,1}-B_{1,1}^{1,0}, \\ \hat D_2 \hat D_3(\varphi(v_1))= -B_{0,1}^{1,1}-B_{0,2}^{0,1}-B_{1,1}^{0,2}, \\
\hat D_1 \hat D_2(\varphi(v_1)) = B_{0,0}^{2,0}-B_{2,0}^{0,2}-B_{1,1}^{0,1}+B_{0,1}^{1,0}, \\ \hat D_2(\varphi(v_1)) = -B_{0,0}^{1,1}-B_{0,1}^{0,1}-B_{1,0}^{0,2}, \\

\hat D_1 \hat D_2 \hat D_3(\varphi(v_1)) = B_{1,1}^{1,1}+B_{0,1}^{2,0}+B_{0,2}^{1,0}, \\ \hat D_3(\varphi(v_1)) = -B_{0,0}^{2,0}+B_{1,1}^{0,1}-B_{1,0}^{1,1}+B_{0,2}^{0,0}.

\end{array}
$$
Поклавши $\varphi(v_1)=0,$ $\varphi(v_2)=0$ знайдемо  необхідні $9$ рівнянь. Роз'язавши в Maple отриману систему із $36$ рівнянь знайдемо значення всіх 27 базисних елементів $B^{i,j}_{k,l}.$
Отже, ми отримали реалізацію незвідного $sl_3$-модуля $B_2^2 \cong \Gamma_{2,2}$  в ${\mathfrak U}(sl_3)a.$ Відповідний  елемент Казиміра буде змішаним конкомітантом
$$
f :=\Delta(B_2^2,S^{2}(X) \cdot S^{2}(U) )=\sum_{\begin{array}{c} \mbox{{\small {\it  i+j}}} \leq 2 \\  \mbox{{\small {\it k+l}}} \leq 2 \end{array} } \, \frac{ 2! 2! B_{i\,j}^{k\,l}x_1^{2-(i+j)} x_2^i x_3^j u_3^{2-(k+l)} u_1^k u_2^l}{i! j! k! l! (2{-}(i+j))! (2-(k+l))!}.
$$
Обчисливши всі коєфіцієнти  $B^{i,j}_{k,l}$ знаходимо явний вигляд $f \colon $
$$
\begin{array}{l}
f :=\Delta(B_2^2,S^{2}(X) \cdot S^{2}(U) )=\\ ({a_{0, \,0}}\,{a_{2, \,0}} - {a_{1, \,0}}^{2})\,x^{2}\,{u_{
3}}^{2} + ({a_{2, \,0}}\,{a_{0, \,2}} - {a_{1, \,1}}^{2})\,x^{2}
\,{u_{1}}^{2} + ({a_{0, \,1}}\,{a_{0, \,3}} - {a_{0, \,2}}^{2})\,
z^{2}\,{u_{2}}^{2} \\
\mbox{} + ({a_{2, \,1}}\,{a_{0, \,3}} - {a_{1, \,2}}^{2})\,z^{2}
\,{u_{1}}^{2} + ({a_{0, \,2}}\,{a_{0, \,0}} - {a_{0, \,1}}^{2})\,
x^{2}\,{u_{2}}^{2} + ({a_{2, \,1}}\,{a_{0, \,1}} - {a_{1, \,1}}^{
2})\,z^{2}\,{u_{3}}^{2} \\
\mbox{} + ( - {a_{2, \,1}}^{2} + {a_{3, \,0}}\,{a_{1, \,2}})\,y^{
2}\,{u_{1}}^{2} + ( - {a_{2, \,0}}^{2} + {a_{1, \,0}}\,{a_{3, \,0
}})\,y^{2}\,{u_{3}}^{2} \\
\mbox{} + ( - {a_{1, \,0}}\,{a_{2, \,0}} + {a_{3, \,0}}\,{a_{0, 
\,0}})\,x\,y\,{u_{3}}^{2} + ({a_{1, \,0}}\,{a_{1, \,2}} - {a_{1, 
\,1}}^{2})\,y^{2}\,{u_{2}}^{2} \\
\mbox{} + ({a_{1, \,2}}\,{a_{0, \,1}} + {a_{1, \,0}}\,{a_{0, \,3}
} - 2\,{a_{1, \,1}}\,{a_{0, \,2}})\,y\,z\,{u_{2}}^{2} \\
\mbox{} + ({a_{1, \,2}}\,{a_{0, \,0}} - 2\,{a_{1, \,1}}\,{a_{0, 
\,1}} + {a_{1, \,0}}\,{a_{0, \,2}})\,x\,y\,{u_{2}}^{2} + (2\,{a_{
1, \,1}}\,{a_{0, \,2}} - 2\,{a_{1, \,2}}\,{a_{0, \,1}})\,{u_{2}}
\,z^{2}\,{u_{3}} \\
\mbox{} + (2\,{a_{1, \,0}}\,{a_{1, \,1}} - 2\,{a_{2, \,0}}\,{a_{0
, \,1}})\,x^{2}\,{u_{3}}\,{u_{1}} + ( - 2\,{a_{1, \,1}}\,{a_{2, 
\,1}} + {a_{2, \,0}}\,{a_{1, \,2}} + {a_{3, \,0}}\,{a_{0, \,2}})
\,x\,y\,{u_{1}}^{2} \\
\mbox{} + (2\,{a_{2, \,0}}\,{a_{2, \,1}} - 2\,{a_{3, \,0}}\,{a_{1
, \,1}})\,y^{2}\,{u_{3}}\,{u_{1}} + (2\,{a_{1, \,1}}\,{a_{1, \,2}
} - 2\,{a_{2, \,1}}\,{a_{0, \,2}})\,z^{2}\,{u_{3}}\,{u_{1}} \\
\mbox{} + ({a_{2, \,0}}\,{a_{0, \,3}} + {a_{2, \,1}}\,{a_{0, \,2}
} - 2\,{a_{1, \,1}}\,{a_{1, \,2}})\,x\,z\,{u_{1}}^{2} + ( - {a_{0
, \,1}}\,{a_{0, \,2}} + {a_{0, \,3}}\,{a_{0, \,0}})\,x\,z\,{u_{2}
}^{2} \\
\mbox{} + (2\,{a_{2, \,0}}\,{a_{1, \,1}} - 2\,{a_{1, \,0}}\,{a_{2
, \,1}})\,{u_{2}}\,y^{2}\,{u_{3}} + ({a_{3, \,0}}\,{a_{0, \,1}}
 + {a_{1, \,0}}\,{a_{2, \,1}} - 2\,{a_{2, \,0}}\,{a_{1, \,1}})\,y
\,z\,{u_{3}}^{2} \\
\mbox{} + (2\,{a_{0, \,1}}\,{a_{1, \,0}} - 2\,{a_{1, \,1}}\,{a_{0
, \,0}})\,{u_{2}}\,x^{2}\,{u_{3}} + (2\,{a_{1, \,1}}\,{a_{2, \,1}
} - 2\,{a_{2, \,0}}\,{a_{1, \,2}})\,y^{2}\,{u_{2}}\,{u_{1}} \\
\mbox{} + ( - 2\,{a_{1, \,0}}\,{a_{0, \,2}} + 2\,{a_{1, \,1}}\,{a
_{0, \,1}})\,x^{2}\,{u_{2}}\,{u_{1}} + ({a_{2, \,1}}\,{a_{0, \,0}
} - 2\,{a_{1, \,0}}\,{a_{1, \,1}} + {a_{2, \,0}}\,{a_{0, \,1}})\,
x\,z\,{u_{3}}^{2} \\
\mbox{} + ( - {a_{2, \,1}}\,{a_{1, \,2}} + {a_{3, \,0}}\,{a_{0, 
\,3}})\,y\,z\,{u_{1}}^{2} + ( - 2\,{a_{1, \,1}}\,{a_{0, \,3}} + 2
\,{a_{1, \,2}}\,{a_{0, \,2}})\,z^{2}\,{u_{2}}\,{u_{1}} \\
\mbox{} + (2\,{a_{1, \,2}}\,{a_{0, \,1}} - 2\,{a_{1, \,0}}\,{a_{0
, \,3}})\,x\,z\,{u_{2}}\,{u_{1}} + ( - 2\,{a_{2, \,1}}\,{a_{0, \,
0}} + 2\,{a_{2, \,0}}\,{a_{0, \,1}})\,{u_{2}}\,x\,y\,{u_{3}} \\
\mbox{} + ( - 2\,{a_{1, \,2}}\,{a_{0, \,0}} + 2\,{a_{1, \,0}}\,{a
_{0, \,2}})\,{u_{2}}\,x\,z\,{u_{3}} \\
\mbox{} + (2\,{a_{1, \,0}}\,{a_{1, \,2}} + 2\,{a_{1, \,1}}^{2} - 
2\,{a_{2, \,1}}\,{a_{0, \,1}} - 2\,{a_{2, \,0}}\,{a_{0, \,2}})\,x
\,z\,{u_{3}}\,{u_{1}} \\
\mbox{} + ( - 2\,{a_{3, \,0}}\,{a_{0, \,2}} + 2\,{a_{2, \,0}}\,{a
_{1, \,2}})\,y\,z\,{u_{3}}\,{u_{1}} \\
\mbox{} + ( - 2\,{a_{2, \,0}}\,{a_{0, \,2}} + 2\,{a_{1, \,1}}^{2}
 - 2\,{a_{1, \,0}}\,{a_{1, \,2}} + 2\,{a_{2, \,1}}\,{a_{0, \,1}})
\,x\,y\,{u_{2}}\,{u_{1}} \\
\mbox{} + (2\,{a_{1, \,0}}\,{a_{2, \,1}} - 2\,{a_{3, \,0}}\,{a_{0
, \,1}})\,x\,y\,{u_{3}}\,{u_{1}} + ( - 2\,{a_{2, \,0}}\,{a_{0, \,
3}} + 2\,{a_{2, \,1}}\,{a_{0, \,2}})\,y\,z\,{u_{2}}\,{u_{1}} \\
\mbox{} + ( - 2\,{a_{1, \,0}}\,{a_{1, \,2}} + 2\,{a_{2, \,0}}\,{a
_{0, \,2}} + 2\,{a_{1, \,1}}^{2} - 2\,{a_{2, \,1}}\,{a_{0, \,1}})
\,{u_{2}}\,y\,z\,{u_{3}}
\end{array}
$$

\end{document}